\documentclass[12pt]{amsart}
\usepackage{amsmath,amssymb}
\newcommand{\IN}{{\mathbb N}}

\newcommand{\cU}{{\mathcal U}}
\DeclareMathOperator{\ad}{ad}
\DeclareMathOperator{\supp}{supp}

\newtheorem{theorem}{Theorem}
\theoremstyle{remark}
\newtheorem{remark}[theorem]{Remark}
\title{Metric spaces with subexponential asymptotic dimension growth}
\author{Narutaka OZAWA}
\address{RIMS, Kyoto University, \mbox{606-8502}}
\email{narutaka@kurims.kyoto-u.ac.jp}
\subjclass{Primary 54F45; Secondary 20F69, 51F99}

\keywords{coarse metric space, asymptotic dimension, Yu's property $A$}
\date{August 08, 2011}

\begin{document}
\begin{abstract}
We prove that a metric space with subexponential asymptotic dimension growth has Yu's property $A$.
\end{abstract}
\maketitle
Let $X$ be a metric space. We denote by $B(x,S)$ the closed ball with center at $x$ and radius $S$.
Recall that $X$ has Yu's property $A$ if there is a net of functions $\zeta^n\colon X\to\ell_1(X)$
such that
\begin{enumerate}
\item $\zeta^n_x\geq0$ and $\|\zeta^n_x\|=1$ for all $n$ and $x$,
\item for every $n$, there is $S_n>0$ such that $\supp\zeta^n_x\subset B(x,S_n)$ for all $x\in X$,
\item\label{c3} for every $R>0$, one has
$\lim_n\sup\{ \|\zeta^n_x-\zeta^n_y\| : d(x,y)\le R\}=0$.
\end{enumerate}
(See Remark~\ref{rem} below on this definition.)
Let us recall from \cite{polygro} the notion of asymptotic dimension growth.
Consider a cover $\cU$ of $X$. It is uniformly bounded if
the diameters of the elements of $\cU$ are uniformly bounded.
The Lebesgue number $L(\cU)$ is the infimum of those $\lambda\geq0$ satisfying
the condition that every subset $B\subset X$ of diameter at most $\lambda$ is contained in some
element of $\cU$. The multiplicity $m(\cU)$ is the maximal number of elements
of $\cU$ with a nonempty intersection.
Finally, the asymptotic dimension function of $X$ is an extended integer valued function $\ad_X$ defined
for $\lambda \in [0,+\infty)$ by
\[
\ad_X(\lambda)=\min\{ m(\cU) :
\cU\mbox{ a uniformly bounded cover of $X$ such that $L(\cU)\geq\lambda$}\}-1.
\]
The metric space $X$ is said to have subexponential asymptotic dimension growth if
$\ad_X$ has subexponential growth. See \cite{polygro} for details.
In this short paper, we prove the following.
\begin{theorem}\label{thm}
A metric space $X$ with subexponential asymptotic dimension growth has Yu's property $A$.
\end{theorem}

This extends Dranishnikov's theorem (\cite{polygro}) that polynomial asymptotic dimension
growth implies Yu's property $A$. In \cite{ds}, Dranishnikov and Sapir introduce
the notion of dimension growth and prove that Thompson's group $F$ has exponential
dimension growth. At the end of \cite{ds}, they ask whether expanders have
exponential dimension growth and whether subexponential dimension growth implies
coarse embeddability. As subexponential dimension growth implies subexponential asymptotic
dimension growth (the converse is also likely, but it is not clear at the moment of writing and in particular
it is not known whether Thompson's group $F$ has exponential asymptotic dimension growth),
the above theorem answers both problems in the affirmative.
(See \cite{willett} for the relationships among Yu's property $A$, coarse embeddability,
and expanders.)

\begin{proof}
The averaging method used in this proof is borrowed from \cite{kaimanovich}.
We will reproduce it here for the reader's convenience.
Let us fix $n$ and construct $\zeta^n\colon X\to\ell_1(X)$.
Take a uniformly bounded cover $\cU=\{ U_i : i\in I\}$ of $X$ such that
$L(\cU)\geq 6n$ and $m(\cU)=\ad_X(6n)+1$.
For every $i\in I$, we choose $x_i\in U_i$ and let
$J\colon\ell_1(I)\to\ell_1(X)$ be the corresponding contraction sending
$\delta_i$ to $\delta_{x_i}$.
For every $x\in X$ and $k\in\IN$, we set
$S_x(k)=\{ i\in I : B(x,k) \subset U_i \}$.
Then, for every $x,y\in X$ and $R\in\IN$ with $d(x,y)\le R$, one has
\[
S_x(k+R) \subset S_x(k) \cap S_y(k) \subset S_x(k) \cup S_y(k) \subset S_x(k-R).
\]
For every non-empty finite subset $S$, we denote by $\xi_S=|S|^{-1}\chi_S$
the uniform probability measure on $S$ and note that
\[
\|\xi_S-\xi_T\|=2\left( 1 - \frac{|S\cap T|}{\max\{|S|,|T|\}}\right).
\]
Now, we define
\[
\eta^n_x=\frac{1}{n}\sum_{k=n+1}^{2n} \xi_{S_x(k)} \in \ell_1(I)
\ \mbox{ and }\
\zeta^n_x=J(\eta^n_x)\in\ell_1(X).
\]
It follows that $\supp\zeta^n_x\subset B(x,\sup_i\mathrm{diam}(U_i))$ for all $x$.
It is left to verify the condition~\ref{c3}.
Let $R\in\IN$ be given and consider $n>R$.
Then, for every $x,y\in X$ with $d(x,y)\le R$, one has
\begin{align*}
\|\eta^n_x-\eta^n_y\|
\le\frac{1}{n}\sum_{k=n+1}^{2n} \| \xi_{S_x(k)} -\xi_{S_y(k)}\|
\le\frac{2}{n}\sum_{k=n+1}^{2n}\left( 1 - \frac{|S_x(k+R)|}{|S_x(k-R)|}\right).
\end{align*}
But since
\begin{align*}
\frac{1}{n}\sum_{k=n+1}^{2n}\frac{|S_x(k+R)|}{|S_x(k-R)|}
\geq \left(\prod_{k=n+1}^{2n}\frac{|S_x(k+R)|}{|S_x(k-R)|}\right)^{1/n}
=\left(\frac{\prod_{k=2n-R+1}^{2n+R}|S_x(k)|}{\prod_{k=n-R+1}^{n+R}|S_x(k)|}\right)^{1/n}
\geq m(\cU)^{-2R/n}
\end{align*}
(here we used the fact that $1\le|S_x(k)|\le m(\cU)$ for all $k\le 2n+R$),
this implies
\[
\|\zeta^n_x-\zeta^n_y\|\le \|\eta^n_x-\eta^n_y\|\le 2\left( 1 - m(\cU)^{-2R/n}\right).
\]
Since $m(\cU)=\ad_X(6n)+1$ has subexponential growth in $n$, we obtain the conclusion.
\end{proof}
\begin{remark}\label{rem}
The definition of Yu's property $A$ used in this paper is different
from the original one (particularly when $X$ does not have bounded geometry).
However, since each $\zeta^n$ constructed in the proof is such that
$m(\cU)!n\zeta^n_x$ are integer-valued for all $x$,
Yu's property $A$ in its original form (\cite{yu}) also follows.
See \cite{willett} for more information.
\end{remark}
\subsection*{Acknowledgment.}
This research was done while the author was visiting at
Texas A{\&}M University in the summer 2011
for ``Workshop in Analysis and Probability'' supported by NSF.
He gratefully acknowledges the kind hospitality and
stimulating environment provided by TAMU and the program organizers.
The author would like to thank Professor M.~Sapir for beautiful
lectures on coarse embeddings during the workshop, in which he asked
if subexponential asymptotic dimension growth implies Yu's property $A$.
The author was partially supported by JSPS and Sumitomo Foundation.

\end{document}